# A Fractional Survival Model

Cheng K. Lee
chengli@uab.edu

Jenq-Daw Lee
Graduate Institute of Political Economy
National Cheng Kung University
Tainan, Taiwan 70101
ROC

**Abstract**
A survival model is derived from the exponential function using the concept of fractional differentiation. The hazard function of the proposed model generates various shapes of curves including increasing, increasing-constant-increasing, increasing-decreasing-increasing, and so-called bathtub hazard curve. The model also contains a parameter that is the maximum of the survival time.

## 1. Introduction
The one-parameter exponential distribution has been applied on many fields. The common form of the probability density function (PDF) is $f(t) = \mu e^{-\mu t}$, where $\mu > 0$ and $t \geq 0$. The cumulative density function (CDF) is $Z(t) = 1 - e^{-\mu t}$. The mean and the variance of this distribution are, respectively, $\frac{1}{\mu}$ and $\frac{1}{\mu^2}$. The hazard function is equal to the first derivative of the CDF divided by the survival function which is equal to 1 minus the CDF also known as the mortality function in survival analysis. Therefore, the hazard function of the one-parameter exponential distribution is $\mu$, a constant. A constant hazard rate does not describe all the observed phenomena in many fields. It is desirable to have a model with a non-constant hazard function. Based on the one-parameter exponential distribution model, we derive a fractional survival model by taking an arbitrary order of the CDF. We then apply the model to 3 cases of real world data.

## 2. The Model

Let *DZ(t)* denote the first derivative of the CDF of the one-parameter exponential function, then the hazard function is $\frac{DZ(t)}{1-Z(t)} = \mu$. Many works have been done to derive models with various hazard rates. Instead of varying the hazard rates, we here employ the idea of Stiassnie (1979) who used a model with an arbitrary order of differentiation to explain the dynamics of viscoelastic materials. Moreover, the arbitrary order is not necessarily an integer. Therefore, we may choose to take the derivative to some arbitrary order of the CDF to be given by $D_t^\lambda Z(t) = \mu e^{-\mu t}$, where $D_t^\lambda Z(t)$ denotes the differentiation of the arbitrary order $\lambda$ with respect to $t$. Using Cauchy formula for repeated integration to solve for $Z(t)$ (see Appendix A1), we have

$Z(t) = \frac{\mu t^\lambda}{\Gamma(\lambda+1)} {}_1F_1[1; \lambda+1; -\mu t]$, where $\Gamma(.)$ is the Gamma function and ${}_1F_1[a;b;c]$ is the confluent hypergeometric function (Gurland, 1958; Muller, 2001) with 3 arguments. However, $Z(t)$ is just the incomplete Gamma distribution, a special case of the confluent hypergeometric function (Luke, 1959). To make the model more general, the first argument of $Z(t)$ is substituted by $a$. Then, $Z(t) = \frac{\mu t^\lambda}{\Gamma(\lambda+1)} {}_1F_1[a; \lambda+1; -\mu t]$. In order for $Z(t)$ to be interpreted as a PDF, it is necessary that, for $t \leq T$, $Z(T) = 1$. After the normalization (see Appendix A2), we obtain the fractional mortality function

$$F(t) = \left(\frac{t}{T}\right)^\lambda \frac{{}_1F_1[a; \lambda+1; -\mu t]}{{}_1F_1[a; \lambda+1; -\mu T]}. \tag{1}$$

The fractional survival function $S(t)$ is, therefore, $1 - F(t)$.

The negative sign in the third argument in equation (2) can be eliminated using Kummer's formula (Gurland, 1958) such that

$$F(t) = e^{\mu(T-t)} \left(\frac{t}{T}\right)^\lambda \frac{{}_1F_1[\lambda+1-a; \lambda+1; \mu t]}{{}_1F_1[\lambda+1-a; \lambda+1; \mu T]}. \tag{2}$$

Taking the first derivative of $F(t)$ in equation (1) and using the differential formula of the confluent hypergeometric function (Abramowitz & Stegun, 1972), the probability density function is

$$f(t) = \frac{\lambda}{T}\left(\frac{t}{T}\right)^{\lambda-1} \frac{{}_1F_1[a; \lambda; -\mu t]}{{}_1F_1[a; \lambda+1; -\mu T]}. \tag{3}$$

Or, using Kummer's formula, the PDF is

$$f(t) = e^{\mu(T-t)} \frac{\lambda}{T} \left(\frac{t}{T}\right)^{\lambda-1} \frac{{}_1F_1[\lambda-a; \lambda; \mu t]}{{}_1F_1[\lambda-a+1; \lambda+1; \mu T]}. \tag{4}$$

The confluent hypergeometric function can be represented by a series and an integral expression (Gurlan, 1958; Muller 2001). Either expression has its own restrictions on the values of the first and the second arguments. In this article, the confluent hypergeometric function is numerically evaluated by Muller's (20001) algorithm based on the series expression in which the second argument cannot be zero or a negative integer. Applying the restrictions to our fractional survival model, $\lambda$ cannot be zero or a negative integer, and $0 \leq t \leq T$.

The hazard function for this fractional survival model provides various shapes of curves (Fig 1) including increasing-constant-increasing (hazard 1 with $\alpha=3$, $\lambda=3$, $\mu=1.5$, $T=20$), decreasing-constant-increasing or so-called bathtub hazard curve (hazard 2 with $\alpha=0.01$, $\lambda=0.01$, $\mu=0.7$, $T=20$), increasing-decreasing-increasing (hazard 3 with $\alpha=3.4$, $\lambda=3.5$, $\mu=2$, $T=20$) and increasing (hazard 4 with $\alpha=-1$, $\lambda=1$, $\mu=2$, $T=20$) hazard rates (Fig 1). However, the mean and the variance of the model do not exist (see Appendix A3).

## 3. Application

Case 1
We use the data of Marriage History File 1958 – 2003 of the Panel Study of Income Dynamics (PSID) to fit the model on time to the first divorce. That is the time in years to the divorce after the first marriage. The data is available at http://simba.isr.umich.edu. In our sample with 12143 individuals born in the 19[th] century, there were 6408 divorces and 5735 marriages still intact when censored in year 2003. The average length of years to divorce of the first marriage was 14.13.
The maximum likelihood estimates and the standard errors based on the second derivatives valued at the maximized log-likelihood function are in Table 1.

Table 1

| Parameter | Estimate | Standard Error |
|---|---|---|
| $\alpha$ | 1.538 | 0.032 |
| $\lambda$ | 1.626 | 0.045 |
| $\mu$ | 0.162 | 0.013 |
| $T$ | 1859.301 | 1932.584 |

The plots of the cumulative hazard curves of the proposed model and the Nelson-Aalen estimates are in figure 2 truncated at time 60. The survival curves of the proposed model and the Nelson-Aalen estimates are in figure 3 truncated at time 60. The increasing-decreasing hazard curve of the fitted model is in figure 4 truncated at time 60. Figure 5 shows the full-scaled hazard curve with increasing-decreasing-increasing rates.

Case 2
We use the same data in case 1 to fit the proposed model on time to first marriage. The sample was censored in 2003 with 17338 married individuals and 6029 unmarried. The average year to the first marriage was 24.01.

The maximum likelihood estimates and the standard errors based on the second derivatives valued at the maximized log-likelihood function are in Table 2.

Table 2

| Parameter | Estimate | Standard Error |
|---|---|---|
| $\alpha$ | 34.025 | 0.551 |
| $\lambda$ | 34.094 | 0.026 |
| $\mu$ | 1.483 | 0.552 |
| $T$ | 94.512 | 0.487 |

The plots of the cumulative hazard curves of the proposed model and the Nelson-Aalen estimates are in figure 6. The survival curves of the proposed model and the Nelson-Aalen estimates are in figure 7. The increasing-decreasing-increasing hazard curve of the fitted model is in figure 8.

Case 3
The data in this application contains survival times in month from 5880 patients after they received coronary artery bypass grafting (CABG). Among these patients, 545 died during the study. The average survival time to death after receiving the procedure was 47 months. The data is available at www.clevelandclinic.org/heartcenter/hazard/default.htm. The maximum likelihood estimates and the standard errors based on the second derivatives valued at the maximized log-likelihood function are in Table 3.

Table 3

| Parameter | Estimate | Standard Error |
|---|---|---|
| $\alpha$ | 0.396 | 0.129 |
| $\lambda$ | 0.302 | 0.023 |
| $\mu$ | -0.021 | 0.003 |
| $T$ | 225.265 | 8.639 |

The plots of the cumulative hazard curves of the proposed model and the Nelson-Aalen estimates are in figure 9. The survival curves of the proposed model and the Nelson-Aalen estimates are in figure 10. The decreasing-constant-increasing hazard curve of the fitted model is in figure 11. The hazard is known as the bathtub curve which is also the three phase hazard function described by Sergeant, Blackstone and Meyns (1997) who analyzed CABG data as well.

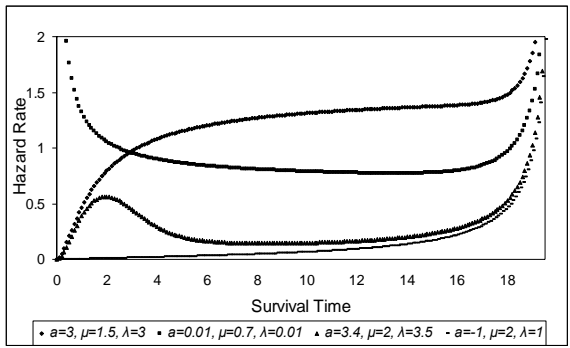

Fig. 1 Hazard curves with various sets of parameters

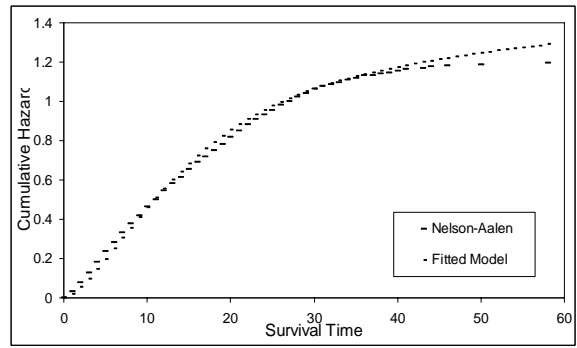

Fig. 2 Cumulative hazards of Nelson-Aalen estimates and the proposed model on time to first divorce truncated at time 60

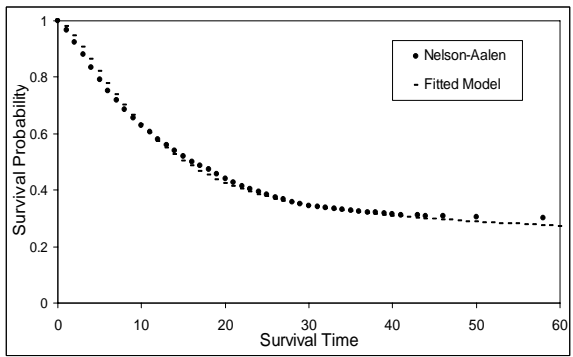

Fig. 3 Survival probabilities of Nelson-Aalen estimates and the proposed model on time to first divorce truncated at time 60

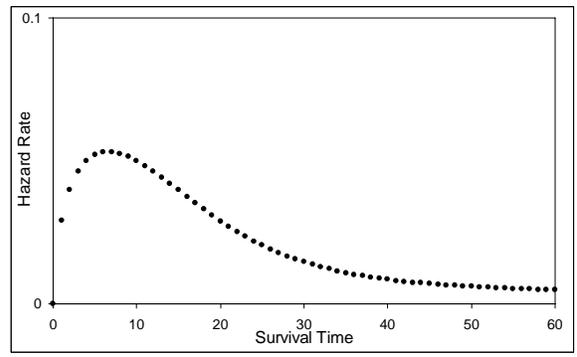

Fig. 4 Hazard rate of the proposed model on time to first divorce truncated at time 60

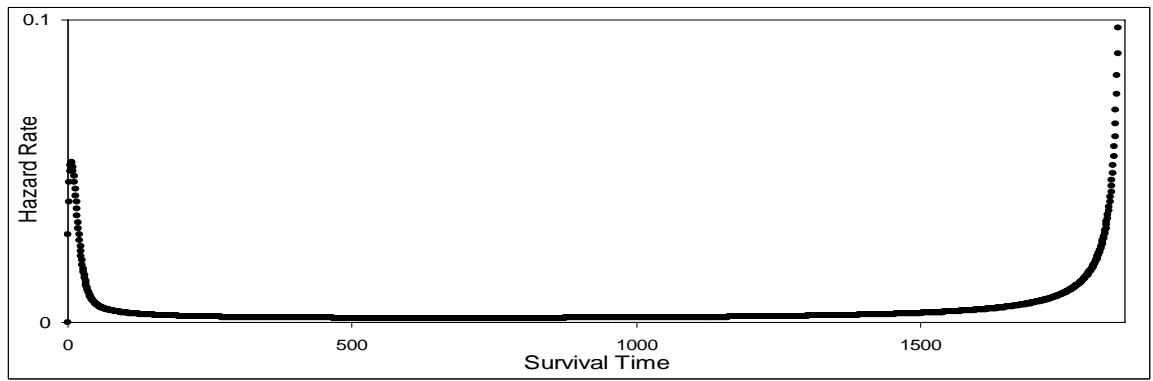

Fig. 5 Full-scaled hazard rate of the proposed model on time to first divorce

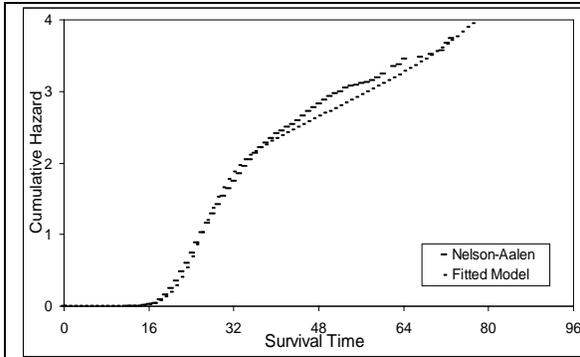

Fig. 6 Cumulative hazards of Nelson-Aalen estimates and the proposed model on time to first marriage

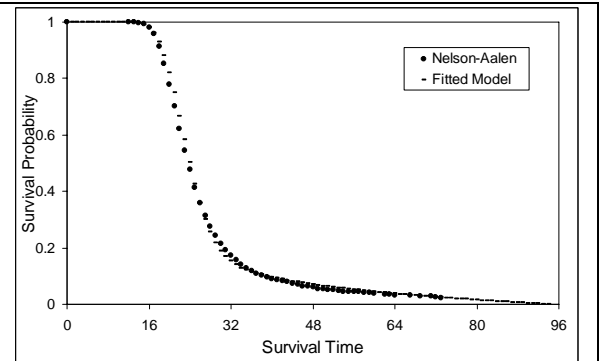

Fig. 7 Survival probabilities of Nelson-Aalen estimates and the proposed model on time to first marriage

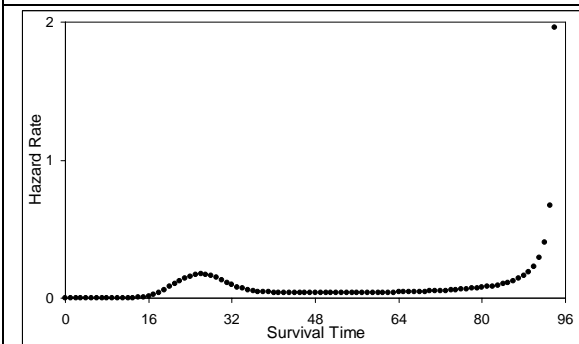

Fig. 8 Hazard rate of the proposed model on time to first marriage

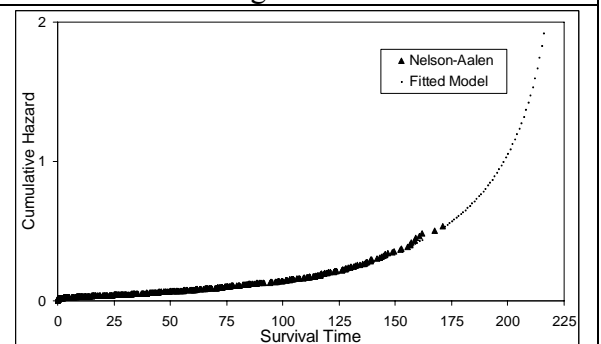

Fig. 9 Cumulative hazards of Nelson-Aalen estimates and the proposed model on survival time after CABG

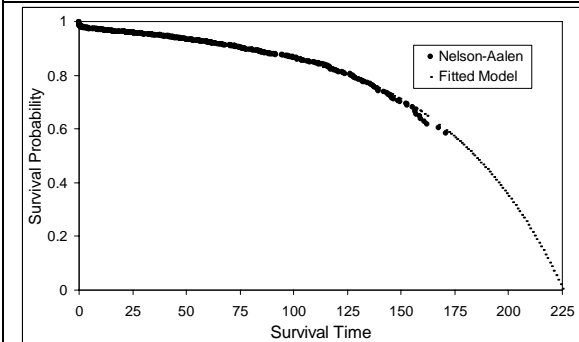

Fig. 10 Survival probabilities of Nelson-Aalen estimates and the proposed model on survival time after CABG

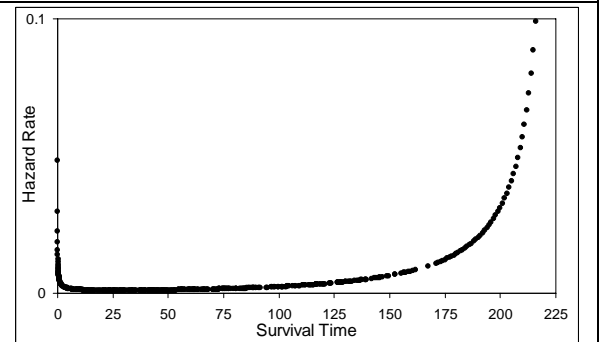

Fig. 11 Hazard rate of the proposed model on survival time after CABG

**Discussion**
In this study, we proposed a survival function with flexible hazard curves with the parameter $T$ that can be regarded as the maximum survival time. However, in case1, the estimated *T* is not significantly different from 0 at 0.05 significance level. It suggests that the model may be overparameterized for the data. A model with an increasing-decreasing hazard curve may be a good candidate. In case 2, the estimated *T* is 94.512 which suggests that the maximum age to the first marriage be 94.512. In case 3, the estimated *T* suggests that the maximum survival month after CABG be 225.265.

**Acknowledgement**

**Appendix**
A1. Solving for $Z(t)$ in Equation 1

$$Z(t) = D_t^{-\lambda}\{D_t^{\lambda}Z(t)\}$$
$$= D_t^{-\lambda}\mu e^{-\mu t}$$
$$= \mu D_t^{-\lambda}e^{-\mu t}$$
$$= \frac{\mu}{\Gamma(\lambda)}\int_0^1 (t-y)^{\lambda-1}e^{-uy}dy, \text{ using Cauchy formula for repeated integration}$$
$$= \frac{\mu}{\Gamma(\lambda)}\int_0^1 (t-ts)^{\lambda-1}e^{-uts}(tds), \text{ let } s=y/t$$
$$= \frac{\mu t^{\lambda}}{\Gamma(\lambda)}\int_0^1 (1-s)^{\lambda-1}e^{-\mu ts}ds$$

The confluent hypergeometric (Gurland, 1958) is defined by series as

$$_1F_1(a;c;z) = \sum_{k=0}^{\infty}\frac{(a)_k z^k}{(c)_k k!}, \text{ where } (a)_k = \frac{\Gamma(a+k)}{\Gamma(a)}.$$

It converges for all real values of *a*, *c*, and *z* and *c* cannot be a negative integer or zero. The confluent hypergeoemtric function can also be presented as an intergral form (Gurland, 1958)

$$_1F_1(a;c;z) = \frac{\Gamma(c)}{\Gamma(a)\Gamma(c-a)}\int_0^1 e^{zt}t^{a-1}(1-t)^{c-a-1}dt, \text{ where } (a)_k = \frac{\Gamma(a+k)}{\Gamma(a)}, \text{ and } c > a > 0.$$

Then,

$$\int_0^1 (1-s)^{\lambda-1}e^{-\mu ts}ds$$
$$= \frac{_1F_1[1;\lambda+1;-\mu t]}{\frac{\Gamma(\lambda+1)}{\Gamma(1)\Gamma(\lambda)}}$$
$$= \frac{_1F_1[1;\lambda+1;-\mu t]}{\lambda}$$

Thus,

$$D_t^{-\lambda}\mu e^{-\mu t} = \frac{\mu t^\lambda}{\lambda \Gamma(\lambda)} {}_1F_1[1; \lambda+1; -\mu t]$$

Therefore, $Z(t) = \frac{\mu t^\lambda}{\Gamma(\lambda+1)} {}_1F_1[1; \lambda+1; -\mu t]$.

A 2. Normalization of $Z(t)$

Let $Z(t) = g(a, \mu, \lambda, t) t^\lambda {}_1F_1[a; \lambda+1; -\mu t]$, where $g(a, \mu, \lambda, t)$ is a function of $a$, $\mu$, $\lambda$, and $t$. The argument $t$ is the normalizing constant. Then, for a constant $T$ ($t \leq T$),
$Z(T) = g(a, \mu, \lambda, T) T^\lambda {}_1F_1[a; \lambda+1; -\mu T]$.
To normalize, for all $T$, $Z(T) = 1 = g(a, \mu, \lambda, T) T^\lambda {}_1F_1[a; \lambda+1; -\mu T]$.

Then, $g(a, \mu, \lambda, T) = \frac{1}{T^\lambda {}_1F_1[a; \lambda+1; -\mu T]}$.

Substituting into $Z(t)$, we obtain

$$Z(t) = \left(\frac{t}{T}\right)^\lambda \frac{{}_1F_1[a; \lambda+1; -\mu t]}{{}_1F_1[a; \lambda+1; -\mu T]}$$

A 3. Moment Generating Function

The moment generating function using the PDF in equation 3 is

$$M_t(s) = \int_0^\infty e^{st} \frac{\lambda}{T}\left(\frac{t}{T}\right)^{\lambda-1} \frac{{}_1F_1[a; \lambda; -\mu t]}{{}_1F_1[a; \lambda+1; -\mu T]} dt$$

$$= \frac{\Gamma(\lambda+1)(-s)^{a-\lambda}(\mu-s)^{-a}}{T^\lambda {}_1F_1[a; \lambda+1; -\mu T]}$$

The first derivative of the moment generating function with respect to $s$ is

$$M_t'(s) = \frac{\Gamma(\lambda+1)(-s)^{a-\lambda-1}(\mu-s)^{-a-1}(\mu\lambda - \mu a - \lambda s)}{T^\lambda {}_1F_1[a; \lambda+1; -\mu T]}$$

The second derivative of the moment generating function with respect to $s$ is
$M_t''(s) =$
$$\frac{\Gamma(\lambda+1)(-s)^{a-\lambda-2}(\mu-s)^{-a-2}\left(s^2\lambda(\lambda+1) + 2s\mu(\lambda+1)(a-\lambda) + (a-\lambda)(a-\lambda-1)\mu^2\right)}{T^\lambda {}_1F_1[a; \lambda+1; -\mu T]}$$

When the confluent hypergeometric is evaluated under the integral expression, $a - \lambda$ must be less than 1 (Gurland, 1958) and, therefore, the first and the second derivative of the moment generating function are undefined. The confluent hypergeometric can also be evaluated under the series expression and both derivatives of the moment generating function are 0. It means that the expectation and the variance are both 0 which implies that $t$ is not a random variable. Therefore, the mean and the variance of the proposed distribution do not exit.